\documentclass{article}
\usepackage{amssymb}
\usepackage{amsmath}

\newcommand{\ie}{i.e.}
\newcommand{\lc}{l.c.}
\newcommand{\nsf}{n.s.f.}

\newcommand{\ot}{\otimes}
\newcommand{\de}{\Delta}
\newcommand{\sde}{\delta}
\newcommand{\Mh}{\hat M}
\newcommand{\deh}{\hat\Delta}

\newcommand{\id}{\text{id}}
\newcommand{\om}{\omega}

\newcommand{\KGA}{\mathcal{K}(G,A)}
\newcommand{\MA}{\text{M}(A)}
\newcommand{\MB}{\text{M}(B)}
\newcommand{\mor}{\text{Mor}(A,B)}

\newcommand{\cs}{C^{*}}

\newcommand{\R}{\mathbb{R}}
\newcommand{\C}{\mathbb{C}}

\newcommand{\Z}{\mathbb{Z}}

\newcommand{\bh}{\mathcal{B}(H)}

\newcommand{\nfs}{n.s.f.}

\def\11{\mbox1\hspace{-.25em}\text{I}}

\newtheorem{thm}{Theorem}
\newtheorem{defi}{Definition}
\newtheorem{lem}{Lemma}
\newtheorem{prop}{Proposition}
\newtheorem{cor}{Corollary}
\newenvironment{pf}{\textbf{Proof.}}{\hfill $\blacksquare$}

\newenvironment{rmk}{\textbf{Remark.}}{}

\begin{document}

\title{A locally compact quantum group of triangular matrices.}
\author{Pierre Fima\footnote{Laboratoire de Math\'ematiques,
Universit\'e de Franche-Comt\'e, 16 route de Gray, 25030 Besancon Cedex, France.
E-mail: fima@math.unicaen.fr} \,\,\,and Leonid Vainerman\footnote{Laboratoire de
Math\'ematiques Nicolas Oresme, Universit\'e de Caen, B.P. 5186, 14032 Caen
Cedex, France. E-mail: vainerman@math.unicaen.fr}}
\date{}
\maketitle

{\bf Dedicated to Professor M.L. Gorbachuk on the occasion of his 70-th anniversary.}

\begin{abstract}
We construct a one parameter deformation of the group of $2\times
2$ upper triangular matrices with determinant $1$ using the
twisting construction. An interesting feature of this new example
of a locally compact quantum group is that the Haar measure is
deformed in a non-trivial way. Also, we give a complete
description of the dual $\cs$-algebra and the dual
comultiplication.
\end{abstract}

\section{Introduction}

In \cite{Enock, Vain}, M. Enock and the second author proposed a
systematic approach to the construction of non-trivial Kac
algebras by twisting. To illustrate it, consider a cocommutative
Kac algebra structure on the group von Neumann algebra
$M=\mathcal{L}(G)$ of a non commutative locally compact (\lc)
group $G$ with comultiplication $\de(\lambda_g)=
\lambda_g\ot\lambda_g$ (here $\lambda_g$ is the left translation
by $g\in G$). Let us define on $M$ another, "twisted",
comultiplication $\de_\Omega(\cdot) = \Omega\de(\cdot)\Omega^*$,
where $\Omega$ is a unitary from $M\ot M$ verifying certain
2-cocycle condition, and construct in this way new, non
cocommutative, Kac algebra structure on $M$. In order to find such
an $\Omega$, let us, following to M. Rieffel \cite{Rief} and M.
Landstad \cite{Land}, take an inclusion $\alpha : L^\infty({\hat
K})\to M$, where $\hat K$ is the dual to some abelian subgroup $K$
of $G$ such that $\delta\vert_{K}=1$, where $\delta(\cdot)$ is the
module of $G$. Then, one lifts a usual 2-cocycle $\Psi$ of $\hat
K:\ \Omega=(\alpha\ot\alpha)\Psi$. The main result of
\cite{Enock}, \cite{Vain} is that the integral by the Haar measure
of $G$ gives also the Haar measure of the deformed object. Recently
P. Kasprzak studied the deformation of l.c. groups by twisting in
\cite{Kas}, and also in this case the Haar measure was not deformed.

In \cite{FimaVain}, the authors extended the twisting construction
in order to cover the case of non-trivial deformation of the
Haar measure. The aim of the present paper is to illustrate this
construction on a concrete example and to compute explicitly all
the ingredients of the twisted quantum group including the dual
$\cs$-algebra and the dual comultiplication. We twist the group von
Neumann algebra $\mathcal{L}(G)$ of the group $G$ of $2\times 2$
upper triangular matrices with determinant $1$ using the abelian
subgroup $K=\C^{*}$ of diagonal matrices of $G$ and a one
parameter family of bicharacters on $K$. In this case, the
subgroup $K$ is not included in the kernel of the modular function
of $G$, this is why the Haar measure is deformed. We compute the
new Haar measure and show that the dual $\cs$-algebra is generated
by $2$ normal operators $\hat{\alpha}$ and $\hat{\beta}$ such that
$$
\hat{\alpha}\hat{\beta}=\hat{\beta}\hat{\alpha}\quad\hat{\alpha}
\hat{\beta}^{*}=q\hat{\beta}^{*}\hat{\alpha},
$$
where $q>0$. Moreover, the comultiplication $\hat{\Delta}$ is given by
$$
\begin{array}{c}
\hat{\Delta}_{t}(\hat{\alpha})=\hat{\alpha}\otimes\hat{\alpha},\
\hat{\Delta}_{t}(\hat{\beta})=\hat{\alpha}\otimes\hat{\beta}
\dot{+}\hat{\beta}\otimes\hat{\alpha}^{-1},
\end{array}
$$
where $\dot{+}$ means the closure of the sum of two operators.

This paper in organized as follows. In Section 2 we recall some
basic definitions and results. In Section 3 we present in detail
our example computing all the ingredients associated. This example
is inspired by \cite{Kas}, but an important difference is that in
the present example the Haar measure is deformed in a non trivial way.
Finally, we collect some useful results in the Appendix.

\section{Preliminaries}

\subsection{Notations}

Let $B(H)$ be the algebra of all bounded linear operators on a
Hilbert space $H$, $\ot$ the tensor product of Hilbert spaces, von
Neumann algebras or minimal tensor product of $\cs$-algebras, and
$\Sigma$ (resp., $\sigma$) the flip map on it. If $H, K$ and $L$
are Hilbert spaces and $X \in B(H \ot L)$ (resp., $X \in B(H \ot
K), X \in B(K \ot L)$), we denote by $X_{13}$ (resp., $X_{12},\
X_{23}$) the operator $(1 \ot \Sigma^*)(X \ot 1)(1 \ot \Sigma)$
(resp., $X\ot 1,\ 1\ot X$) defined on $H \ot K \ot L$. For any
subset $X$ of a Banach space $E$, we denote by $\langle X\rangle$
the vector space generated by $X$ and $[X]$ the closed vector
space generated by $X$. All l.c. groups considered in this paper
are supposed to be second countable, all Hilbert spaces are
separable and all von Neumann algebras have separable preduals.

Given a {\it normal semi-finite faithful} (n.s.f.) weight $\theta$
on a von Neumann algebra $M$ (see \cite{S}), we denote: ${\mathcal M}^+_\theta
 = \{ x \in M^+ \mid \theta(x) < + \infty \},\ {\mathcal N}_\theta =
\{ x \in M \mid x^*x \in M^+_\theta \},$ and ${\mathcal
M}_\theta=\langle {\mathcal M}^+_\theta\rangle$.

When $A$ and $B$ are $\cs$-algebras, we denote by M$(A)$ the algebra
of the multipliers of $A$ and by Mor$(A,B)$ the set of the morphisms from
$A$ to $B$.

\subsection{$G$-products and their deformation}

For the notions of an action of a l.c. group $G$ on a $\cs$-algebra $A$,
a $\cs$ dynamical system $(A,G,\alpha)$, a crossed product $G\,_{\alpha}
\ltimes A$ of $A$ by $G$ see \cite{ped}. The crossed product has the
following universal property:

For any $\cs$-covariant representation $(\pi, u, B)$ of $(A,G,\alpha)$ (here
$B$ is a $\cs$-algebra, $\pi:A\to B$ a morphism, $u$ is a group morphism from
$G$ to the unitaries of $M(B)$, continuous for the strict topology), there is a
unique morphism $\rho\in\text{Mor}(G\,_{\alpha}\ltimes A,B)$ such that
$$\rho(\lambda_{t})=u_{t},\quad\rho(\pi_{\alpha}(x))=\pi(x)\quad\forall t\in
G,x\in A.$$

\begin{defi}
Let $G$ be a l.c. abelian group, $B$ a $\cs$-algebra, $\lambda$ a morphism from $G$ to the unitary group of $\MB$, continuous in the strict topology of $\MB$, and $\theta$ a continuous action of $\hat{G}$ on $B$. The triplet $(B,\lambda,\theta)$ is called a $G$-product if $\theta_{\gamma}(\lambda_{g})= \overline{\langle\gamma,g\rangle}\lambda_{g}$ for all $\gamma\in\hat{G}$, $g\in G$.
\end{defi}
The unitary representation $\lambda\,:\, G\rightarrow\MB$ generates a morphism : $$\lambda\in\text{Mor}(\cs(G),B).$$
Identifying $\cs(G)$ with $C_{0}(\hat{G})$, one gets a morphism $\lambda\in\text{Mor}(C_{0}(\hat{G}),B)$ which is
defined in a unique way by its values on the characters
$$
u_{g}=(\gamma\mapsto \langle\gamma,g\rangle)\in C_{b}(\hat{G}):\
\lambda(u_{g})=\lambda_{g},\quad\text{for all}\,\,g\in G.
$$
One can check that $\lambda$ is injective.

The action $\theta$ is done by: $\theta_{\gamma}(\lambda(u_{g}))=\theta_{\gamma}(\lambda_{g})=
\overline{\langle\gamma,g\rangle}\lambda_{g}
=\lambda(u_{g}(.-\gamma)).$
Since the $u_{g}$ generate $C_{b}(\hat{G})$, one deduces that:
$$
\theta_{\gamma}(\lambda(f))=\lambda(f(.-\gamma)),\quad\text{for all}\,\,f\in C_{b}(\hat{G}).
$$

The following definition is equivalent to the original definition by Landstad \cite{Land} (see \cite{Kas}):

\begin{defi}\label{defland}
Let $(B,\lambda,\theta)$ be a $G$-poduct and $x\in\MB$. One says that $x$ verifies the Landstad conditions if
\begin{equation}
\left\{\begin{array}{ll}
(i) & \theta_{\gamma}(x)=x,\quad\text{for any}\quad \gamma\in\hat{G},\\\label{condland}
(ii) & \text{the application}\,\, g\mapsto\lambda_{g}x\lambda_{g}^{*}\,\,\text{is continuous},\\
(iii) & \lambda(f)x\lambda(g)\in B,\quad\text{for any}\quad f,g\in C_{0}(\hat{G}).\\
\end{array}\right.
\end{equation}
\end{defi}
The set $A\in\MB$ verifying these conditions is a $\cs$-algebra called \textit{the Landstad algebra} of the $G$-product
$(B,\lambda,\theta)$. Definition \ref{defland} implies that if $a\in A$, then $\lambda_{g}a\lambda_{g}^{*}\in A$ and
the map $g\mapsto\lambda_{g}a\lambda_{g}^{*}$ is continuous. One gets then an action of $G$ on $A$.

One can show that the inclusion $A\to\MB$ is a morphism of $\cs$-algebras, so
$\MA$ can be also included into $\MB$. If $x\in\MB$, then $x\in\MA$ if and
only if
\begin{equation}
\left\{
\begin{array}{ll}
(i) & \theta_{\gamma}(x)=x,\,\,\text{for
all}\,\,\gamma\in\hat{G},\\\label{multland}
(ii) & \text{for all}\,\,a\in A,\,\,\text{the application}\,\, g\mapsto\lambda_{g}x\lambda_{g}^{*}a \,\,\text{is continuous.}\\
\end{array}
\right.
\end{equation}
Let us note that two first conditions of $(\ref{condland})$ imply $(\ref{multland})$.

The notions of $G$-product and crossed product are closely
related. Indeed, if $(A,G,\alpha)$ is a $\cs$-dynamical system
with $G$ abelian, let $B=G\,_{\alpha}\ltimes A$ be the crossed
product and $\lambda$ the canonical morphism from $G$ into the
unitary group of $\MB$, continuous in the strict topology, and
$\pi\in\mor$ the canonical morphism of $C^*$-algebras. For
$f\in\KGA$ and $\gamma\in\hat{G}$, one defines
$(\theta_{\gamma}f)(t)=\overline{\langle\gamma,t\rangle}f(t)$. One
shows that $\theta_{\gamma}$ can be extended to the automorphisms
of $B$ in such a way that $(B,\hat{G},\theta)$ would be a
$\cs$-dynamical system. Moreover, $(B,\lambda,\theta)$ is a
$G$-product and the associated Landstad algebra is $\pi(A)$.
$\theta$ is called \textit{the dual action}. Conversely, if
$(B,\lambda,\theta)$ is a $G$-product, then one shows that there
exists a $\cs$-dynamical system $(A,G,\alpha)$ such that
$B=G\,_{\alpha}\ltimes A$. It is unique (up to a covariant
isomorphism), $A$ is
 the Landstad algebra of $(B,\lambda,\theta)$ and $\alpha$ is the action of $G$ on $A$ given by $\alpha_{t}(x)=\lambda_{t}x\lambda_{t}^{*}$.

\begin{lem}\cite{Kas}\label{lem36Kas}
Let $(B,\lambda,\theta)$ be a $G$-product and $V\subset A$ be a vector subspace of the Landstad algebra such that:
\begin{itemize}
\item $\lambda_{g}V\lambda_{g}^{*}\subset V$, for any $g\in G$,
\item $\lambda(C_{0}(\hat{G}))V\lambda(C_{0}(\hat{G}))$ is dense in $B$.
\end{itemize}
Then $V$ is dense in $A$.
\end{lem}

Let $(B,\lambda,\theta)$ be a $G$-product, $A$ its Landstad algebra, and $\Psi$ a continuous bicharacter on $\hat{G}$. For $\gamma\in\hat{G}$, the function on $\hat{G}$ defined by $\Psi_{\gamma}(\omega)= \Psi(\omega,\gamma)$ generates a family of unitaries $\lambda(\Psi_{\gamma})\in\MB$. The bicharacter condition implies:
$$
\theta_{\gamma}(U_{\gamma_{2}})=
\lambda(\Psi_{\gamma_{2}}(.-\gamma_{1}))=
\overline{\Psi(\gamma_{1},\gamma_{2})}
U_{\gamma_{2}},\quad\forall\gamma_{1},\gamma_{2}\in\hat{G}.
$$
One gets then a new action $\theta^{\Psi}$ of $\hat{G}$ on $B$:
$$
\theta^{\Psi}_{\gamma}(x)=U_{\gamma}\theta(x)U_{\gamma}^{*}.
$$
Note that, by commutativity of $G$, one has:
$$\theta^{\Psi}_{\gamma}(\lambda_{g})=
U_{\gamma}\theta(\lambda_{g})U_{\gamma}^{*}=
\overline{\langle\gamma,g\rangle}\lambda_{g},\quad\forall \gamma\in\hat{G},g\in G.
$$
The triplet $(B,\lambda,\theta^{\Psi})$ is then a $G$-product, called a \textit{deformed} $G$-product.

\subsection{Locally compact quantum groups \cite{VaesLC}, \cite{KV2}}

A pair $(M,\de)$ is called a (von Neumann algebraic) l.c.\ quantum
group  when
\begin{itemize}
\item $M$ is a von Neumann algebra and $\de : M \to M \ot M$ is a
normal and unital $*$-homomorphism which is coassociative: $(\de
\ot \id)\de = (\id \ot \de)\de$ (i.e., $(M,\de)$ is a Hopf-von
Neumann algebra). \item There exist n.s.f. weights $\varphi$ and
$\psi$ on $M$ such that
\begin{itemize}
\item $\varphi$ is left invariant in the sense that $\varphi
\bigl( (\om \ot \id)\de(x) \bigr) = \varphi(x) \om(1)$ for all $x
\in {\mathcal M}_{\varphi}^+$ and $\om \in M_*^+$, \item $\psi$ is
right invariant in the sense that $\psi \bigl( (\id \ot \om)\de(x)
\bigr) = \psi(x) \om(1)$ for all $x \in {\mathcal M}_{\psi}^+$ and
$\om \in M_*^+$.
\end{itemize}
\end{itemize}
Left and  right invariant weights are unique up to a positive
scalar.

Let us represent $M$ on the GNS Hilbert space of $\varphi$ and
define a unitary $W$ on $H \ot H$ by
$$
W^* (\Lambda(a) \ot \Lambda(b)) = (\Lambda \ot \Lambda)(\de(b)(a
\ot 1)), \quad\text{for all}\; a,b \in N_{\phi}\; .
$$
Here, $\Lambda$ denotes the canonical GNS-map for $\varphi$,
$\Lambda \ot \Lambda$ the similar map for $\varphi \ot \varphi$.
One proves that $W$ satisfies the {\it pentagonal equation}:
$W_{12} W_{13} W_{23} = W_{23} W_{12}$, and we say that $W$ is a
{\it multiplicative unitary}. The von Neumann algebra $M$ and the
comultiplication on it can be given in terms of $W$ respectively
as
$$
M = \{ (\id \ot \om)(W) \mid \om \in B(H)_* \}^{-\sigma-strong*} \;
$$
and $\de(x) = W^* (1 \ot x) W$, for all $x \in M$. Next, the l.c.\
quantum group $(M,\de)$ has an antipode $S$, which is the unique
$\sigma$-strongly* closed linear map from $M$ to $M$ satisfying
$(\id \ot \om)(W) \in {\mathcal D}(S)$ for all $\om \in B(H)_*$
and $S(\id \ot \om)(W) = (\id \ot \om)(W^*)$ and such that the
elements $(\id \ot \om)(W)$ form a $\sigma$-strong* core for $S$.
$S$ has a polar decomposition $S = R \tau_{-i/2}$, where $R$ (the
unitary antipode) is an anti-automorphism of $M$ and $\tau_t$ (the
scaling group of $(M,\de)$) is a strongly continuous
one-parameter group of automorphisms of $M$. We
have $\sigma (R \ot R) \de = \de R$, so $\varphi R$ is a right
invariant weight on $(M,\de)$ and we take $\psi:= \varphi R$.

Let $\sigma_t$ be the modular automorphism group of $\varphi$.
There exist a number $\nu > 0$, called the scaling constant, such
that $\psi \, \sigma_t = \nu^{-t} \, \psi$ for all $t \in
\mathbb{R}$. Hence (see \cite{VaesRN}), there is a unique positive, self-adjoint
operator $\sde_M$ affiliated to $M$, such that $\sigma_t(\sde_M) =
\nu^t \, \sde_M$ for all $t \in \mathbb{R}$ and $\psi =
\varphi_{\sde_M}$. It is called the modular element of $(M,\de)$.
If $\sde_M=1$ we call $(M,\de)$ unimodular. The scaling constant
can be characterized as well by the relative invariance $\varphi
\, \tau_t = \nu^{-t} \, \varphi$.

For the dual l.c.\ quantum group $(\Mh,\deh)$ we have :
$$\Mh = \{(\om \ot \id)(W) \mid \om \in B(H)_*
\}^{-\sigma-strong*}$$ and $\deh(x) = \Sigma W (x \ot 1) W^*
\Sigma$ for all $x \in \Mh$. A left invariant
\nsf{} weight $\hat \varphi$ on $\Mh$ can be constructed
explicitly and the associated multiplicative unitary is $\hat{W} =
\Sigma W^* \Sigma$.

Since $(\Mh,\deh)$ is again a l.c.\ quantum group, let us denote its
antipode by $\hat{S}$, its unitary antipode by $\hat{R}$ and its
scaling group by $\hat{\tau}_t$. Then we can construct the dual of
$(\Mh,\deh)$, starting from the left invariant weight
$\hat\varphi$. The bidual l.c.\ quantum group $(\hat\Mh,\hat\deh)$
is isomorphic to $(M,\de)$.

$M$ is commutative if and only if $(M,\de)$ is generated by a usual
l.c.\ group $G : M=L^{\infty}(G), (\de_G f)(g,h) = f(gh),$ $(S_Gf)(g)
= f(g^{-1}),\ \varphi_G(f)=\int f(g)\; dg$, where $f\in L^{\infty}(G),
\ g,h\in G$ and we integrate with respect to the left Haar measure
$dg$ on $G$. Then $\psi_G$ is given by $\psi_G(f) = \int f(g^{-1}) \;
dg$ and $\sde_M$ by the strictly positive function $g \mapsto
\sde_G(g)^{-1}$.

$L^\infty(G)$ acts on $H=L^2(G)$ by multiplication and
$(W_G\xi)(g,h)=\xi(g,g^{-1}h),$ for all $\xi\in H\ot H=L^2(G\times
G)$. Then $\Mh=\mathcal L(G)$ is the group von Neumann algebra
generated by the left translations $(\lambda_g)_{g\in G}$ of $G$
and $\deh_G(\lambda_g)= \lambda_g\ot\lambda_g$. Clearly,
$\deh_G^{op}:=\sigma\circ\deh_G=\deh_G$, so $\deh_G$ is
cocommutative.

$(M,\Delta)$ is a Kac algebra (see \cite{ES}) if $\tau_{t}=\id$, for
all $t$, and $\delta_{M}$ is affiliated with the center of $M$. In
particular, this is the case when $M=L^{\infty}(G)$ or $M=
\mathcal{L}(G)$.

We can also define the $\cs$-algebra of continuous functions vanishing
at infinity on $(M,\Delta)$ by
$$
A=\left[(\id\otimes\omega)(W)\,\,|\,\,\omega\in\bh_{*}\right]
$$
and the reduced $\cs$-algebra (or dual $\cs$-algebra) of $(M,\Delta)$ by
$$
\hat{A}=\left[(\omega\otimes\id)(W)\,\,|\,\,\omega\in\bh_{*}\right].
$$
In the group case we have $A=C_{0}(G)$ and $\hat{A}=C_{r}(G)$. Moreover,
we have $\Delta\in\text{Mor}(A,A\otimes A)$ and
$\hat{\Delta}\in\text{Mor}(\hat{A},\hat{A}\otimes\hat{A})$.

A l.c. quantum group is called compact if $\varphi(1_M)<\infty$ and
discrete if its dual is compact.

\subsection{Twisting of locally compact quantum groups
\cite{FimaVain}}

Let $(M,\de)$ be a locally compact quantum group and $\Omega$ a unitary
in $M\ot M$. We say that $\Omega$ is a 2-cocycle on $(M,\de)$ if
$$(\Omega\ot 1)(\de\ot \id)(\Omega)=(1\ot\Omega)(\id\ot\de)(\Omega).$$
As an example we can consider $M=L^\infty(G)$, where $G$ is a l.c.
group, with $\de_G$ as above, and $\Omega=\Psi(\cdot,\cdot)\in
L^\infty(G\times G)$ a usual 2-cocycle on $G$, i.e., a mesurable
function with values in the unit circle $\mathbb T\subset \mathbb
C$ verifying $$\Psi(s_1,s_2)\Psi(s_1 s_2,s_3)=
\Psi(s_2,s_3)\Psi(s_1,s_2 s_3),\ \text{for almost all}\,\,
s_1,s_2,s_3\in G.$$ This is the case for any measurable
bicharacter on $G$.

When $\Omega$ is a $2$-cocycle on $(M,\Delta)$, one can check that
$\de_\Omega (\cdot)= \Omega\de(\cdot)\Omega^*$ is a new
coassociative comultiplication on $M$. If $(M,\de)$ is discrete
and $\Omega$ is any 2-cocycle on it, then $(M,\de_\Omega)$ is
again a l.c. quantum group (see \cite{BichonVaes},
finite-dimensional case was treated in \cite{Vain}). In the
general case, one can proceed as follows. Let
$\alpha:(L^\infty(G),\de_G) \to (M,\de)$ be an inclusion of
Hopf-von Neumann algebras, i.e., a faithful unital normal
*-homomorphism such that $(\alpha\ot\alpha)\circ\de_G=
\de\circ\alpha$. Such an inclusion allows to construct a 2-cocycle
of $(M,\de)$ by lifting a usual 2-cocycle of $G :
\Omega=(\alpha\ot\alpha)\Psi$. It is shown in \cite{Enock} that if
the image of $\alpha$ is included into the centralizer of the left
invariant weight $\varphi$, then $\varphi$ is also left invariant
for the new comultiplication $\Delta_{\Omega}$.

In particular, let $G$ be a non commutative \lc{} group and $K$ a
closed abelian subgroup of $G$. By Theorem 6 of \cite{TT}, there exists
a faithful unital normal *-homomorphism
$\hat{\alpha}\,:\,\mathcal{L}(K)\rightarrow\mathcal{L}(G)$ such that
$$
\hat{\alpha}(\lambda^{K}_{g})=\lambda_{g},\quad\text{for all}\,\,g\in
K,\quad\text{and}\quad\hat\Delta\circ\hat\alpha=(\hat\alpha\otimes
\hat\alpha)\circ\hat\Delta_{K},
$$
where $\lambda^{K}$ and $\lambda$ are the left regular representation
of $K$ and $G$ respectively, and $\hat\Delta_K$ and $\hat\Delta$ are the
comultiplications on $\mathcal{L}(K)$ and $\mathcal{L}(G)$ repectively.
The composition of $\hat{\alpha}$ with the canonical isomorphism
$L^{\infty}(\hat{K})\simeq \mathcal{L}(K)$ given by the Fourier
tranformation, is a faithful unital normal *-homomorphism
$\alpha\,:\,L^{\infty}(\hat{K})\rightarrow \mathcal{L}(G)$ such that
$\Delta\circ\alpha=(\alpha\otimes\alpha)\circ\Delta_{\hat{K}},$ where
$\Delta_{\hat{K}}$ is the comultiplication on $L^{\infty}(\hat{K})$. The
left invariant weight on $\mathcal{L}(G)$ is the Plancherel weight for
which
$$
\sigma_{t}(x)=\delta_{G}^{it}x\delta_{G}^{-it},\quad\text{for all}\,\,x\in
\mathcal{L}(G),
$$ where $\delta_{G}$ is the modular function of $G$. Thus, $\sigma_{t}(\lambda_{g})=\delta_{G}^{it}(g)\lambda_{g}$ or
$$
\sigma_{t}\circ\alpha(u_{g})=\alpha(u_{g}(\cdot-\gamma_{t})),
$$
where $u_{g}(\gamma)=\langle\gamma,g\rangle,\ g\in G,\gamma\in\hat G,\
\gamma_{t}$ is the character $K$ defined by $\langle\gamma_{t},g\rangle=
\delta_{G}^{-it}(g).$ By linearity and density we obtain:
$$
\sigma_{t}\circ\alpha(F)=\alpha(F(\cdot-\gamma_{t})),\quad\text{for all}
\,\,F\in L^{\infty}(\hat{K}).
$$
This is why we do the following assumptions. Let $(M,\Delta)$ be a l.c.
quantum group, $G$ an abelian \lc{} group and $\alpha:(L^\infty(G),\de_G)
\to (M,\de)$ an inclusion of Hopf-von Neumann algebras. Let $\varphi$ be
the left invariant weight, $\sigma_{t}$ its modular group, $S$ the antipode,
$R$ the unitary antipode, $\tau_{t}$ the scaling group. Let $\psi=\varphi
\circ R$ be the right invariant weight and $\sigma_{t}^{'}$ its modular
group. Also we denote by $\delta$ the modular element of $(M,\Delta)$.
Suppose that there exists a continuous group homomorphism $t\mapsto
\gamma_{t}$ from $\R$ to $G$ such that
$$
\sigma_{t}\circ\alpha(F)=\alpha(F(\cdot-\gamma_{t})),
\quad\text{for all}\,\,F\in L^{\infty}(G).
$$
Let $\Psi$ be a continuous bicharacter on $G$. Notice that
$(t,s)\mapsto\Psi(\gamma_{t},\gamma_{s})$ is a continuous
bicharacter on $\mathbb{R}$, so there exists $\lambda
>0$ such that $\Psi(\gamma_{t},\gamma_{s})=\lambda^{ist}$. We
define:
$$
u_{t}=\lambda^{i\frac{t^{2}}{2}}\alpha\left(\Psi(.,-\gamma_{t})\right)\quad
\text{and}\quad v_{t}=\lambda^{i\frac{t^{2}}{2}}\alpha\left(\Psi(-\gamma_{t},.)
\right).
$$
The $2$-cocycle equation implies that $u_{t}$ is a $\sigma_t$-cocyle and $v_t$
is a $\sigma^{'}_{t}$-cocycle. The Connes' Theorem gives two \nfs{}
weights on $M$, $\varphi_{\Omega}$ and $\psi_{\Omega}$, such that
$$
u_{t}=[D\varphi_{\Omega}\,:\,D\varphi]_{t}\quad\text{and}\quad
v_{t}=[D\psi_{\Omega}\,:\,D\psi]_{t}.
$$
The main result of \cite{FimaVain} is as follows:

\begin{thm}
$(M,\Delta_{\Omega})$ is a \lc{} quantum group with left and right
invariant weight $\varphi_{\Omega}$ and $\psi_{\Omega}$
respectively. Moreover, denoting by a subscript or a superscript
$\Omega$ the objects associated with $(M,\Delta_{\Omega})$ one
has:
\begin{itemize}
\item $\tau_{t}^{\Omega}=\tau_{t}$,
\item $\nu_{\Omega}=\nu$ and
$\delta_{\Omega}=\delta A^{-1}B$,
\item
$\mathcal{D}(S_{\Omega})=\mathcal{D}(S)$ and, for all
$x\in\mathcal{D}(S),\ S_{\Omega}(x)=uS(x)u^{*}.$
\end{itemize}
\end{thm}

Remark that, because $\Psi$ is a bicharacter on $G$,
$t\mapsto\alpha(\Psi(.,-\gamma_{t}))$ is a representation of $\R$
in the unitary group of $M$ and there exists a positive
self-adjoint operator $A$ affiliated with $M$ such that
$$\alpha(\Psi(.,-\gamma_{t}))=A^{it},\quad\text{for all}\,\,t\in\R.$$
We can also define a positive self-adjoint operator $B$ affiliated with $M$ such that
$$\alpha(\Psi(-\gamma_{t},.))=B^{it}.$$
We obtain :
$$
u_{t}=\lambda^{i\frac{t^{2}}{2}}A^{it},\quad v_{t}=\lambda^{i\frac{t^{2}}{2}}B^{it}.
$$
Thus, we have $\varphi_{\Omega}=\varphi_{A}$ and $\psi_{\Omega}=\psi_{B}$, where $\varphi_{A}$ and $\psi_{B}$ are the weights defined by S. Vaes in \cite{VaesRN}.

One can also compute the dual $\cs$-algebra and the dual comultiplication. We put:
$$
L_{\gamma}=\alpha(u_{\gamma}),\quad R_{\gamma}=JL_{\gamma}J,\quad\text{for all}\,\,\gamma\in\hat{G}.
$$
From the representation $\gamma\mapsto L_{\gamma}$ we get the
unital *-homomorphism $\lambda_{L}\,:\,L^{\infty}(G)\rightarrow M$
and from the representation $\gamma\mapsto R_{\gamma}$ we get the
unital normal *-homomorphism
$\lambda_{R}\,:\,L^{\infty}(G)\rightarrow M^{'}$. Let $\hat{A}$ be
the reduced $\cs$-algebra of $(M,\Delta)$. We can define an action
of $\hat{G}^{2}$ on $\hat{A}$ by
$$
\alpha_{\gamma_{1},\gamma_{2}}(x)=
L_{\gamma_{1}}R_{\gamma_{2}}xR_{\gamma_{2}}^{*}L_{\gamma_{1}}^{*}.
$$
Let us consider the crossed product $\cs$-algebra $B=\hat{G}^{2}\,_{\alpha} \ltimes
\hat{A}.$ We will denote by $\lambda$ the canonical morphism from $\hat{G}^{2}$ to
the unitary group of $M(B)$ continuous in the strict topology on M$(B)$,
$\pi\in$Mor$(\hat{A},B)$ the canonical morphism and $\theta$ the dual action of
$G^{2}$ on $B$. Recall that the triplet $(\hat{G}^{2},\lambda,\theta)$ is a
$\hat{G}^{2}$-product. Let us denote by $(\hat{G}^{2},\lambda,\theta^{\Psi})$ the
$\hat{G}^{2}$-product obtained by deformation of the $\hat{G}^{2}$-product
$(\hat{G}^{2},\lambda,\theta)$ by the bicharacter $\omega(g,h,s,t):=
\overline{\Psi(g,s)}\Psi(h,t)$ on $G^{2}$.

The dual deformed action $\theta^{\Psi}$ is done by
$$
\theta^{\Psi}_{(g_{1},g_{2})}(x)=U_{g_{1}}V_{g_{2}}
\theta_{(g_{1},g_{2})}(x)U_{g_{1}}^{*}V_{g_{2}}^{*},\quad
\text{for any}\,\,g_{1},g_{2}\in G,\,\,x\in B,
$$
where $U_{g}=\lambda_{L}(\Psi_{g}^{*})$, $V_{g}=\lambda_{R}(\Psi_{g})$, $\Psi_{g}(h)=\Psi(h,g)$.

Considering $\Psi_{g}$ as an element of $\hat{G}$, we get a morphism from $G$ to $\hat{G}$, also
noted $\Psi$, such that $\Psi(g)=\Psi_{g}$. With these notations, one has  $U_{g}=u_{(\Psi(-g),0)}$ and $V_{g}=u_{(0,\Psi(g))}$. Then the action $\theta^{\Psi}$ on $\pi(\hat{A})$ is done by
\begin{equation}\label{DualActHatM}
\theta^{\Psi}_{(g_{1},g_{2})}(\pi(x))=\pi(\alpha_{(\Psi(-g_{1}),\Psi(g_{2}))}(x)).
\end{equation}

Let us consider the Landstad algebra $A^{\Psi}$ associated with this $\hat{G}^{2}$-product.
By definition of $\alpha$ and the universality of the crossed product we get a morphism
\begin{equation}\label{EqDualSurj}
\rho\in\text{Mor}(B,\mathcal{K}(H)),\quad\rho(\lambda_{\gamma_{1},\gamma_{2}})=
L_{\gamma_{1}}R_{\gamma_{2}}\quad\text{et}\quad\rho(\pi(x))=x.
\end{equation}
It is shown in \cite{FimaVain} that $\rho(A^{\Psi})=\hat{A}_{\Omega}$ and that $\rho$ is injective on $A^{\Psi}$.
This gives a canonical isomorphism $A^{\Psi}\simeq \hat{A}_{\Omega}$. In the sequel we identify
$A^{\Psi}$ with $\hat{A}_{\Omega}$. The comultiplication can be described in the following way.
First, one can show that, using universality of the crossed product, there exists a unique morphism
$\Gamma\in\text{Mor}(B,B\otimes B)$ such that:
$$
\Gamma\circ\pi=(\pi\otimes\pi)\circ\hat{\Delta}\quad
\text{and}\quad\Gamma(\lambda_{\gamma_{1},\gamma_{2}})=
\lambda_{\gamma_{1},0}\otimes \lambda_{0,\gamma_{2}}.
$$
Then we introduce the unitary
$\Upsilon=(\lambda_{R}\otimes\lambda_{L})(\tilde{\Psi})\in
\text{M}(B\otimes B)$, where $\tilde{\Psi}(g,h)=\Psi(g,gh)$. This
allows us to define the *-morphism
$\Gamma_{\Omega}(x)=\Upsilon\Gamma(x)\Upsilon^{*}$ from $B$ to
M$(B\otimes B)$. One can show that
$\Gamma_{\Omega}\in\text{Mor}(A^{\Psi},A^{\Psi}\otimes A^{\Psi})$
is the comultiplication on $A^{\Psi}$.

Note that if $M=\mathcal{L}(G)$ and $K$ is an abelian closed subgroup of $G$,
the action $\alpha$ of $K^{2}$ on $C_{0}(G)$ is the left-right action.


\section{Twisting of the group of $2\times 2$ upper triangular
matrices with determinant $1$}


Consider the following subgroup of $SL_{2}(\mathbb{C})$ :
$$
G:=\left\{\left(
\begin{array}{cc}
z & \omega \\
0 & z^{-1} \\
\end{array}\right)\,\, ,\,\,\, z\in\mathbb{C}^{*}\, ,\,\,\omega\in\mathbb{C}\right\}.
$$
Let $K\subset G$ be the subgroup of diagonal matrices in $G$, \ie{} $K=\C^{*}$. The elements of $G$ will be denoted by
$(z,\omega)$, $z\in\mathbb{C}$, $\omega\in\mathbb{C}^{*}$.
The modular function of $G$ is
$$
\delta_{G}((z,\omega))=|z|^{-2}.
$$
Thus, the morphism $(t\mapsto\gamma_{t})$ from $\R$ to $\widehat{\C^{*}}$ is given by
$$
\langle \gamma_{t}, z\rangle=|z|^{2it},\quad\text{for all}\,\,z\in\C^{*},t\in\R.
$$
We can identify $\widehat{\mathbb{C}^{*}}$ with $\mathbb{Z}\times\mathbb{R}_{+}^{*}$ in the following way:
$$
\mathbb{Z}\times\mathbb{R}_{+}^{*}\rightarrow\widehat{\mathbb{C}^{*}},\quad (n,\rho)\mapsto\gamma_{n,\rho}=(re^{i\theta}\mapsto e^{i\ln{r}\ln{\rho}}e^{in\theta}).
$$
Under this identification, $\gamma_{t}$ is the element $(0,e^{t})$ of $\mathbb{Z}\times\mathbb{R}_{+}^{*}$. For all $x\in\R$, we define a bicharacter on $\mathbb{Z}\times\mathbb{R}_{+}^{*}$ by
$$
\Psi_{x}((n,\rho),(k,r))=e^{ix(k\ln{\rho}-n\ln{r})}.
$$
We denote by $(M_{x},\Delta_{x})$ the twisted \lc{} quantum group. We have:
$$
\Psi_{x}((n,\rho),\gamma_{t}^{-1})=e^{ixtn}=u_{e^{ixt}}((n,\rho)).
$$
In this way we obtain the operator $A_{x}$ deforming the Plancherel weight:
$$
A_{x}^{it}=\alpha(u_{e^{ixt}})=\lambda^{G}_{(e^{itx},0)}.
$$
In the same way we compute the operator $B_{x}$ deforming the Plancherel weight:
$$
B_{x}^{it}=\lambda^{G}_{(e^{-ixt},0)}=A_{x}^{-it}.
$$
Thus, we obtain for the modular element :
$$
\delta_{x}^{it}=A_{x}^{-it}B_{x}^{it}=\lambda^{G}_{(e^{-2itx},0)}.
$$
The antipode is not deformed. The scaling group is trivial but, if
$x\neq 0$, $(M_{x},\Delta_{x})$ is not a Kac algebra because
$\delta_{x}$ is not affiliated with the center of $M$. Let us look
if $(M_{x},\Delta_{x})$ can be isomorphic for different values of
$x$. One can remark that, since $\Psi_{-x}=\Psi_{x}^{*}$ is
antisymmetric and $\Delta$ is cocommutative, we have
$\Delta_{-x}=\sigma\Delta_{x},$ where $\sigma$ is the flip on
$\mathcal{L}(G)\otimes\mathcal{L}(G)$. Thus,
$(M_{-x},\Delta_{-x})\simeq (M_{x},\Delta_{x})^{\text{op}},$ where
"op" means the opposite quantum group. So, it suffices to treat
only strictly positive values of $x$. The twisting deforms only
the  comultiplication, the weights and the modular element. The
simplest invariant distinguishing the $(M_{x},\Delta_{x})$ is then
the specter of the modular element. Using the Fourier
transformation in the first variable, on has immediately
$\text{Sp}(\delta_{x})=q_{x}^{\Z}\cup\{0\},$ where
$q_{x}=e^{-2x}$. Thus, if $x\neq y,\ x>0,y>0$, one has
$q_{x}^{\Z}\neq q_{y}^{\Z}$ and, consequently,
$(M_{x},\Delta_{x})$ and $(M_{y},\Delta_{y})$ are non isomorphic.

We compute now the dual $\cs$-algebra. The action of $K^{2}$
on $C_{0}(G)$ can be lifted to its Lie algebra $\mathbb{C}^{2}$. The lifting
does not change the result of the deformation (see \cite{Kas}, Proposition
3.17) but simplify calculations. The action of $\mathbb{C}^{2}$ on $C_{0}(G)$
will be denoted by $\rho$. One has
\begin{equation}\label{gd1}
\rho_{z_{1},z_{2}}(f)(z,\omega)=f(e^{z_{2}-z_{1}}z,e^{-(z_{1}+z_{2})}\omega).
\end{equation}
The group $\mathbb{C}$ is self-dual, the duality is given by
$$
(z_{1},z_{2})\mapsto \exp\left(i\text{Im}(z_{1}z_{2})\right).
$$
The generators $u_{z}$, $z\in\C$, of $C_{0}(\C)$ are given by
$$
u_{z}(w)=\exp\left(i\text{Im}(zw)\right),\quad z,w\in\C.
$$
Let $x\in\mathbb{R}$. We will consider the following bicharacter on
$\mathbb{C}$:
$$
\Psi_{x}(z_{1},z_{2})=\exp\left(ix\text{Im}(z_{1}\overline{z}_{2})\right).
$$
Let $B$ be the crossed product $\cs$-algebra $\mathbb{C}^{2}\ltimes
C_{0}(G)$. We denote by
$((z_{1},z_{2})\mapsto\lambda_{z_{1},z_{2}})$ the canonical group
homomorphism from $G$ to the unitary group of $\MB$, continuous
for the strict topology, and $\pi\in$Mor$(C_{0}(G),B)$ the canonical
homomorphism. Also we denote by
$\lambda\in\text{Mor}(C_{0}(G^{2}),B)$ the morphism given by the
representation $((z_{1},z_{2})\mapsto\lambda_{z_{1},z_{2}})$. Let
$\theta$ be the dual action of $\mathbb{C}^{2}$ on $B$. We have,
for all $z,w\in\C$, $\Psi_{x}(w,z)=u_{x\overline{z}}(w)$. The
deformed dual action is given by
\begin{equation}\label{actiondualedeformee}
\theta_{z_{1},z_{2}}^{\Psi_{x}}(b)=
\lambda_{-x\overline{z}_{1},x\overline{z}_{2}}\theta_{z_{1},z_{2}}(b)
\lambda^{*}_{-x\overline{z}_{1},x\overline{z}_{2}}.
\end{equation}
Recall that
\begin{equation}\label{dualelambda}
\theta_{z_{1},z_{2}}^{\Psi_{x}}(\lambda(f))=
\theta_{z_{1},z_{2}}(\lambda(f))=\lambda(f(\cdot-z_{1},\cdot-
z_{2})),\quad\forall f\in C_{b}(\C^{2}).
\end{equation}
Let $\hat{A}_{x}$ be the associated Landstad algebra. We identify
$\hat{A}_{x}$ with the reduced $\cs$-algebra of
$(M_{x},\Delta_{x})$. We will now construct two normal operators
affiliated with $\hat{A}_{x}$, which generate $\hat{A}_{x}$. Let
$a$ and $b$ be the coordinate functions on $G$, and
$\alpha=\pi(a)$, $\beta=\pi(b)$. Then $\alpha$ and $\beta$ are
normal operators, affiliated with $B$, and one can see,
using $(\ref{gd1})$, that
\begin{equation}\label{lambdaalphabeta}
\lambda_{z_{1},z_{2}}\alpha\lambda_{z_{1},z_{2}}^{*}=e^{z_{2}-z_{1}}\alpha,\quad
\lambda_{z_{1},z_{2}}\beta\lambda_{z_{1},z_{2}}^{*}=e^{-(z_{1}+z_{2})}\beta.
\end{equation}
We can deduce, using $(\ref{actiondualedeformee})$, that
\begin{equation}\label{actionalbe}
\theta^{\Psi_{x}}_{z_{1},z_{2}}(\alpha)=e^{x(\overline{z}_{1}+\overline{z}_{2})}\alpha\,\,
,\,\,\,\,\theta^{\Psi_{x}}_{z_{1},z_{2}}(\beta)=e^{x(\overline{z}_{1}-\overline{z}_{2})}\beta.
\end{equation}
Let $T_{l}$ and $T_{r}$ be the infinitesimal generators of the
left and right shift respectively, \ie{} $T_{l}$ and $T_{r}$ are
normal, affiliated with $B$, and
$$
\lambda_{z_{1},z_{2}}= \exp\left(i\text{Im}(z_{1}T_{l})\right)
\exp\left(i\text{Im}(z_{2}T_{r})\right),\quad\text{for
all }\,\,z_{1},z_{2}\in\C.
$$
Thus, we have:
$$
\lambda(f)=f(T_{l},T_{r}),\quad\text{for all}\,\,f\in
C_{b}(\mathbb{C}^{2}).
$$
Let $U=\lambda(\Psi_{x})$, we define the
following normal operators affiliated with $B$:
$$
\hat{\alpha}=U^{*}\alpha U\,\, ,\,\,\,\, \hat{\beta}=U\beta U^{*}.
$$

\begin{prop}\label{generators1}
The operators $\hat{\alpha}$ and $\hat{\beta}$ are affiliated with
$\hat{A}_{x}$ and generate $\hat{A}_{x}$.
\end{prop}

\begin{pf}
First let us show that $f(\hat{\alpha}), f(\hat{\beta})\in
\text{M}(\hat{A}_{t})$, for all $f\in C_{0}(\mathbb{C})$. One has,
using $(\ref{dualelambda})$:
\begin{eqnarray*}
\theta^{\Psi_{x}}_{z_{1},z_{2}}(U) &=& \lambda\left(\Psi_{x}(.-z_{1},.-z_{2})\right)\\
&=& U e^{ix\text{Im}(-\overline{z}_{2}T_{l})}
e^{ix\text{Im}(\overline{z}_{1}T_{r})}\Psi_{x}(z_{1},z_{2})\\
&=&U
\lambda_{-x\overline{z}_{2},x\overline{z}_{1}}\Psi_{x}(z_{1},z_{2}).
\end{eqnarray*}
Now, using (\ref{actionalbe}) and (\ref{lambdaalphabeta}), we
obtain:
$$
\theta^{\Psi_{x}}_{z_{1},z_{2}}(\hat{\alpha})=\hat{\alpha}
,\quad \theta^{\Psi_{x}}_{z_{1},z_{2}}(\hat{\beta})=\hat{\beta},
\quad\text{for all}\,\, z_{1},z_{2}\in\C.
$$

Thus, for all $f\in C_{0}(\C)$, $f(\hat{\alpha})$ and
$f(\hat{\beta})$ are fixed points for the action $\theta^{\Psi_{x}}$.
Let $f\in C_{0}(\mathbb{C})$. Using $(\ref{lambdaalphabeta})$ we
find:
\begin{eqnarray}\notag
\lambda_{z_{1},z_{2}}f(\hat{\alpha})\lambda_{z_{1},z_{2}}^{*}
&=&U^{*}f(e^{z_{2}-z_{1}}\alpha)U,\\
\lambda_{z_{1},z_{2}}f(\hat{\beta})\lambda_{z_{1},z_{2}}^{*}
&=&U^{*}f(e^{-(z_{1}+z_{2})}\beta)U.\label{eqc*2}
\end{eqnarray} Because $f$ is continuous and vanish at infinity,
the applications
$$
(z_{1},z_{2})\mapsto\lambda_{z_{1},z_{2}}
f(\hat{\alpha})\lambda_{z_{1},z_{2}}^{*}\quad\text{and}\quad
(z_{1},z_{2})\mapsto\lambda_{z_{1},z_{2}}f(\hat{\beta})\lambda_{z_{1},z_{2}}^{*}
$$
are norm-continuous and $f(\hat{\alpha}),f(\hat{\beta})\in
M(\hat{A}_{x})$, for all $f\in C_{0}(\C)$.

Taking in mind Proposition \ref{propaff} (see Appendix), in order to show
that $\hat{\alpha}$ is affiliated with $\hat{A}_{x}$, it suffices to show that
the vector space $\mathcal{I}$
generated by $f(\hat{\alpha})a$, with $f\in C_{0}(\mathbb{C})$ and $a \in
\hat{A}_{x}$, is dense in $\hat{A}_{x}$. Using $(\ref{eqc*2})$, we see that
$\mathcal{I}$ is globally invariant under the action implemented by $\lambda$.
Let $g(z)=(1+\overline{z}z)^{-1}$. As $\lambda(C_{0}(\mathbb{C}^{2}))U=\lambda(C_{0}(\mathbb{C}^{2})),$
we can deduce that the closure of
$\lambda(C_{0}(\mathbb{C}^{2}))g(\hat{\alpha})\hat{A}_{x}
\lambda(C_{0}(\mathbb{C}^{2}))$
is equal to
$$
\left[\lambda(C_{0}(\mathbb{C}^{2}))(1+\alpha^{*}\alpha)^{-1}U^{*}\hat{A}_{x}
\lambda(C_{0}(\mathbb{C}^{2}))\right].
$$
As the set $U^{*}\hat{A}_{x}\lambda(C_{0}(\mathbb{C}^{2}))$ is dense in $B$ and
$\alpha$ is affiliated with $B$, the set $\lambda(C_{0}(\mathbb{C}^{2}))(1+\alpha^{*}\alpha)^{-1}U^{*}\hat{A}_{x}
\lambda(C_{0}(\mathbb{C}^{2}))$ is dense in $B$. Moreover, it is included in
$\lambda(C_{0}(\mathbb{C}^{2}))\mathcal{I}\lambda(C_{0}(\mathbb{C}^{2}))$,
so $\lambda(C_{0}(\mathbb{C}^{2}))\mathcal{I}\lambda(C_{0}(\mathbb{C}^{2}))$
is dense in $B$. We conclude, using Lemma \ref{lem36Kas}, that
$\mathcal{I}$ is dense in $\hat{A}_{x}$. One can show in the same
way that $\hat{\beta}$ is affiliated with $\hat{A}_{x}$.

Now, let us show that $\hat{\alpha}$ and $\hat{\beta}$ generate
$\hat{A}_{x}$. By Proposition \ref{propgen}, it suffices to
show that
$$
\mathcal{V}=\left\langle f(\hat{\alpha})g(\hat{\beta}),\,\,f,g\in
C_{0}(\mathbb{C})\right\rangle
$$ is a dense vector subspace of $\hat{A}_{x}$. We have shown above that
the elements of $\mathcal{V}$ satisfy the two first Landstad's conditions. Let
$$
\mathcal{W}=\left[\lambda(C_{0}(\mathbb{C}^{2}))
\mathcal{V}\lambda(C_{0}(\mathbb{C}^{2}))\,\right].
$$
We will show that $\mathcal{W}=B$. This proves that the elements of
$\mathcal{V}$ satisfy the third Landstad's condition, and then
$\mathcal{V}\subset\hat{A}_{x}$. Then $(\ref{eqc*2})$ shows that $\mathcal{V}$
is globally invariant under the action implemented by $\lambda$, so $\mathcal{V}$
is dense in $\hat{A}_{x}$ by Lemma \ref{lem36Kas}. One has:
$$
\mathcal{W}=\left[xU^{*}f(\alpha)U^{2}g(\beta)U^{*}y\, ,\,\,f,g\in C_{0}(\mathbb{C}),\, x,y\in \lambda(C_{0}(\mathbb{C}^{2}))\right].
$$
Because $U$ is unitary, we can substitute $x$ with $xU$ and $y$ with $Uy$
without changing $\mathcal{W}$:
$$\mathcal{W}=\left[xf(\alpha)U^{2}g(\beta)y\, ,\,\,f,g\in C_{0}(\mathbb{C}),\, x,y\in \lambda(C_{0}(\mathbb{C}^{2}))\right].$$
Using, for all $f\in C_{0}(\C)$, the norm-continuity of the application
$$
(z_{1},z_{2})\mapsto\lambda_{z_{1},z_{2}}f(\alpha)
\lambda_{z_{1},z_{2}}^{*}=e^{z_{2}-z_{1}}\alpha,
$$
one deduces that
\begin{eqnarray*}
&&\left[f(\alpha)x\, ,\,\,f\in C_{0}(\mathbb{C}),\,x\in \lambda(C_{0}(\mathbb{C}^{2}))\right]\\
&&=\left[xf(\alpha)\, ,\,\,f\in C_{0}(\mathbb{C}),\,x\in
\lambda(C_{0}(\mathbb{C}^{2}))\right].
\end{eqnarray*}
In particular,
$$
\mathcal{W}=\left[f(\alpha)xU^{2}g(\beta)y\,
,\,\,f,g\in C_{0}(\mathbb{C}),\, x,y\in
\lambda(C_{0}(\mathbb{C}^{2}))\right].
$$
Now we can commute $g(\beta)$ and $y$, and we obtain:
$$
\mathcal{W}=\left[f(\alpha)xU^{2}yg(\beta)\, ,\,\,f,g\in C_{0}(\mathbb{C}),\, x,y\in \lambda(C_{0}(\mathbb{C}^{2}))\right].
$$
Substituting $x\mapsto xU^{*}$, $y\mapsto U^{*}y$, one has:
$$
\mathcal{W}=\left[f(\alpha)xyg(\beta)\, ,\,\,f,g\in C_{0}(\mathbb{C}),\, x,y\in \lambda(C_{0}(\mathbb{C}^{2}))\right].
$$
Commuting back $f(\alpha)$ with $x$ and $g(\beta)$ with $y$, we obtain:
$$
\mathcal{W}=\left[xf(\alpha)g(\beta)y\, ,\,\,f,g\in C_{0}(\mathbb{C}),\, x,y\in \lambda(C_{0}(\mathbb{C}^{2}))\right]=B.
$$
This concludes the proof.
\end{pf}

We will now find the commutation relations between $\hat{\alpha}$
and $\hat{\beta}$.
\begin{prop}\label{diagrel}
One has:
\begin{enumerate}
\item $\alpha$ et $T^{*}_{l}+T^{*}_{r}$ strongly commute and $\hat{\alpha}=e^{x(T_{l}^{*}+T_{r}^{*})}\alpha$.
\item $\beta$ et $T^{*}_{l}-T^{*}_{r}$ strongly commute and $\hat{\beta}=e^{x(T_{l}^{*}-T_{r}^{*})}\beta$.
\end{enumerate}
Thus, the polar decompositions are given by :
\begin{eqnarray*}
&\text{Ph}(\hat{\alpha})=e^{-ix\text{Im}(T_{l}+T_{r})}\text{Ph}(\alpha)\,
,\quad
|\hat{\alpha}|=e^{x\text{Re}(T_{l}+T_{r})}|\alpha|,&\\
&\text{Ph}(\hat{\beta})=e^{-ix\text{Im}(T_{l}-T_{r})}\text{Ph}(\beta)\,
,\quad |\hat{\beta}|=e^{x\text{Re}(T_{l}-T_{r})}|\beta|.&
\end{eqnarray*}
Moreover, we have the following relations:
\begin{enumerate}
\item $|\hat{\alpha}|$ and $|\hat{\beta}|$ strongly commute,
\item $\text{Ph}(\hat{\alpha})\text{Ph}(\hat{\beta})=
    \text{Ph}(\hat{\beta})\text{Ph}(\hat{\alpha})$,
\item $\text{Ph}(\hat{\alpha})|\hat{\beta}|\text{Ph}(\hat{\alpha})^{*}=e^{4x}|\hat{\beta}|$,
\item $\text{Ph}(\hat{\beta})|\hat{\alpha}|\text{Ph}(\hat{\beta})^{*}=e^{4x}|\hat{\alpha}|$.
\end{enumerate}
\end{prop}

\begin{pf}
Using $(\ref{lambdaalphabeta})$, we find, for all $z\in\C$:
$$
e^{i\text{Im}\left(z(T_{l}^{*}+T_{r}^{*})\right)}\alpha e^{-i\text{Im}\left(z(T_{l}^{*}+T_{r}^{*})\right)}=
\lambda_{-\overline{z},-\overline{z}}\alpha\lambda_{-\overline{z},-\overline{z}}^{*}
=e^{-\overline{z}+\overline{z}}\alpha=\alpha.
$$ Thus,
$T_{l}^{*}+T_{r}^{*}$ and $\alpha$ strongly commute. Moreover,
because $e^{ix\text{Im}T_{l}T_{l}^{*}}=1$, one has:
$$
\hat{\alpha}=e^{-ix\text{Im}T_{l}T_{r}^{*}}\alpha e^{ix\text{Im}T_{l}T_{r}^{*}}
=e^{-ix\text{Im}T_{l}(T_{l}+T_{r})^{*}}\alpha
e^{ix\text{Im}T_{l}(T_{l}+T_{r})^{*}}.
$$
We can now prove the point $1$ using the equality
$e^{-ix\text{Im}T_{l}\omega}\alpha
e^{ix\text{Im}T_{l}\omega}=e^{x\omega}\alpha$, the preceding
equation and the fact that $T_{l}^{*}+T_{r}^{*}$ and $\alpha$
strongly commute. The proof of the second assertion is similar and
the polar decompositions follows. From $(\ref{lambdaalphabeta})$
we deduce :

\begin{eqnarray*}
e^{-ix\text{Im}(T_{r}-T_{l})}\alpha e^{ix\text{Im}(T_{r}-T_{l})}
&=&e^{-2x}\alpha,\\
e^{ix\text{Im}(T_{l}+T_{r})}\beta e^{-ix\text{Im}(T_{l}+T_{r})} &=& e^{-2x}\beta,\\
e^{ix\text{Re}(T_{r}-T_{l})}\alpha e^{-ix\text{Re}(T_{r}-T_{l})}&=& e^{2ix}\alpha,\\
e^{ix\text{Re}(T_{l}+T_{r})}\beta e^{-ix\text{Re}(T_{l}+T_{r})}
&=& e^{-2ix}\beta.
\end{eqnarray*}
It is now easy to prove the last relations from the preceding
equations and the polar decompositions.
\end{pf}

We can now give a formula for the comultiplication.

\begin{prop}\label{diagcom}
Let $\hat{\Delta}_{x}$ be the comultiplication on $\hat{A}_{x}$. One has:
$$
\begin{array}{c}
\hat{\Delta}_{x}(\hat{\alpha})=\hat{\alpha}\otimes\hat{\alpha},\
\hat{\Delta}_{x}(\hat{\beta})=
\hat{\alpha}\otimes\hat{\beta}\dot{+}\hat{\beta}\otimes\hat{\alpha}^{-1}.
\end{array}
$$
\end{prop}

\begin{pf}
Using the Preliminaries, we have that
$\hat{\Delta}_{x}=\Upsilon\Gamma(.)\Upsilon^{*}$, where
$$
\Upsilon=e^{ix\text{Im}T_{r}\otimes T_{l}^{*}}
$$ and $\Gamma$ is given by
\begin{itemize}
\item $\Gamma(T_{l})=T_{l}\otimes 1$,  $\Gamma(T_{r})=1\otimes
T_{r}$; \item $\Gamma$ restricted to $C_{0}(G)$ is equal to the
comultiplication $\Delta_G$.
\end{itemize}
Define $R=\Upsilon\Gamma(U^{*})$. One has
$\Delta_{x}(\hat{\alpha})=R(\alpha\otimes\alpha)R^{*}$. Thus, it
is sufficient to show that $(U\otimes U)R$ commute with
$\alpha\otimes\alpha$. Indeed, in this case, one has
$$\hat{\Delta}_{x}(\hat{\alpha})=R(\alpha\otimes\alpha)R^{*}
=(U^{*}\otimes U^{*})(U\otimes
U)R(\alpha\otimes\alpha)R^{*}(U^{*}\otimes U^{*})(U\otimes U)
=\hat{\alpha}\otimes\hat{\alpha}.$$ Let us show that $(U\otimes
U)R$ commute with $\alpha\otimes\alpha$. From the equality
$U=e^{ix\text{Im}T_{l}T_{r}^{*}}$, we deduce that
$$
\Gamma(U^{*})=e^{-ix\text{Im}T_{l}\otimes T_{r}^{*}},\quad U\otimes U=e^{ix\text{Im}(T_{l}T_{r}^{*}\otimes 1+1\otimes T_{l}T_{r}^{*})}.
$$
Thus, $R=e^{-ix\text{Im}(T_{r}^{*}\otimes T_{l}+T_{l}\otimes
T_{r}^{*})}$ and
$$(U\otimes U)R=e^{ix\text{Im}(T_{l}T_{r}^{*}\otimes 1+1\otimes T_{l}T_{r}^{*}-T_{r}^{*}\otimes T_{l}-T_{l}\otimes T_{r}^{*})}.$$
Notice that
$$
T_{l}T_{r}^{*}\otimes 1+1\otimes T_{l}T_{r}^{*}-T_{r}^{*}
\otimes T_{l}-T_{l}\otimes T_{r}^{*}
=(T_{l}\otimes 1-1\otimes T_{l})(T_{r}^{*}\otimes 1-1\otimes
T_{r}^{*}).
$$
Thus, it suffices to show that $T_{l}\otimes 1-1\otimes T_{l}$ and
$T_{r}^{*}\otimes 1-1\otimes T_{r}^{*}$ strongly commute with
$\alpha\otimes\alpha$. This follows from the equations
\begin{eqnarray*}
&&e^{i\text{Im}z(T_{r}^{*}\otimes 1-1\otimes
T_{r}^{*})}(\alpha\otimes\alpha)e^{-i\text{Im}z(T_{r}^{*}\otimes
1-1\otimes T_{r}^{*})}\\
&&=(\lambda_{0,-\overline{z}}\otimes\lambda_{0,\overline{z}})
(\alpha\otimes\alpha)(\lambda_{0,-\overline{z}}\otimes\lambda_{0,\overline{z}})^{*}\\
&&=e^{-\overline{z}}e^{\overline{z}}\alpha\otimes\alpha=\alpha\otimes\alpha,\quad\forall z\in\C\\
\end{eqnarray*}
and
\begin{eqnarray*}
&&e^{i\text{Im}z(T_{l}\otimes 1-1\otimes
T_{l})}(\alpha\otimes\alpha)e^{-i\text{Im}z(T_{l}\otimes
1-1\otimes
T_{l})}\\
&&=(\lambda_{z,0}\otimes\lambda_{-z,0})(\alpha\otimes\alpha)
(\lambda_{z,0}\otimes\lambda_{-z,0})^{*}\\
&&=e^{-z}e^{z}\alpha\otimes\alpha=\alpha\otimes\alpha,\quad\forall z\in\C.\\
\end{eqnarray*}
Put $S=\Upsilon\Gamma(U)$. One has:
$$
\hat{\Delta}_{x}(\hat{\beta})=S(\alpha\otimes\beta+\beta\otimes\alpha^{-1})S^{*}
=S(\alpha\otimes\beta)S^{*}\dot{+}S(\beta\otimes\alpha^{-1})S^{*}.
$$
As before, we see that it suffices to show that $(U\otimes
U^{*})S$ commutes with $\alpha\otimes\beta$ and that $(U^{*}\otimes
U)S$ commutes with $\beta\otimes\alpha^{-1}$, and one can check this
in the same way .
\end{pf}

Let us summarize the preceding results in the following corollary (see \cite{Wor3, Kas}
for the definition of commutation relation between unbounded operators):

\begin{cor}
Let $q=e^{8x}$. The $\cs$-algebra $\hat{A}_{x}$ is generated by
$2$ normal operators $\hat{\alpha}$ and $\hat{\beta}$ affiliated
with $\hat{A}_{x}$ such that
$$
\hat{\alpha}\hat{\beta}=
\hat{\beta}\hat{\alpha}\quad\hat{\alpha}\hat{\beta}^{*}=q\hat{\beta}^{*}\hat{\alpha}.
$$
Moreover, the comultiplication $\hat{\Delta}_{x}$ is given by
$$
\begin{array}{c}
\hat{\Delta}_{x}(\hat{\alpha})=\hat{\alpha}\otimes\hat{\alpha},\
\hat{\Delta}_{x}(\hat{\beta})=\hat{\alpha}\otimes\hat{\beta}\dot{+}\hat{\beta}\otimes\hat{\alpha}^{-1}.
\end{array}
$$
\end{cor}

\begin{rmk}
One can show, using the results of \cite{FimaVain}, that the application $(q\mapsto W_{q})$ which maps the parameter $q$ to the multiplicative unitary of the twisted l.c. quantum group is continuous in the $\sigma$-weak topology.
\end{rmk}

\section{Appendix}


Let us cite some results on operators affiliated with a $\cs$-algebra.

\begin{prop}\label{propaff}
Let $A\subset\bh$ be a non degenerated $\cs$-subalgebra and $T$ a
normal densely defined closed operator on $H$. Let $\mathcal{I}$
be the vector space generated by $f(T)a$, where $f\in C_{0}(\C)$
and $a\in A$. Then:
$$\left(T\eta A\right) \Leftrightarrow
\left(\begin{array}{c}
f(T)\in\MA\,\,\text{for any}\,\,f\in C_{0}(\C)\\
\text{et}\,\,\mathcal{I}\,\,\text{is dense in}\,\,A\\
\end{array}\right).$$
\end{prop}

\begin{pf}
If $T$ is affiliated with $A$, then it is clear that
$f(T)\in\text{M}(A)$ for any $f\in C_{0}(\C)$, and that
$\mathcal{I}$ is dense in $A$ (because $\mathcal{I}$ contains
$(1+T^{*}T)^{-\frac{1}{2}}A$). To show the converse, consider the
*-homomorphism $\pi_{T}: C_{0}(\C)\rightarrow\MA$ given by
$\pi_{T}(f)=f(T)$. By hypothesis, $\pi_{T}(C_{0}(\C))A$ is dense
in $A$. So, $\pi_{T}\in\text{Mor}(C_{0}(\C),A)$ and
$T=\pi_{T}(z\mapsto z)$ is then affiliated with $A$.
\end{pf}

\begin{prop}\label{propgen}
Let $A\subset\bh$ be a non degenerated $\cs$-subalgebra and
$T_{1}, T_{2},\ldots, T_{N}$
 normal operators affiliated with
$A$. Let us denote by $\mathcal{V}$ the vector space generated by
the products of the form $f_{1}(T_{1})f_{2}(T_{2})\ldots
f_{N}(T_{N})$, with $f_{i}\in C_{0}(\C)$. If $\mathcal{V}$ is a
dense vector subspace of $A$, then $A$ is generated by $T_{1},
T_{2},\ldots, T_{N}$.
\end{prop}

\begin{pf}
This follows from Theorem 3.3 in \cite{Wor2}.
\end{pf}

\end{document}